\documentclass[12pt,reqno]{amsart}
\usepackage{bm}
\theoremstyle{plain} 
\newtheorem{thm}{Theorem}[section]

\newtheorem{lem}[thm]{Lemma}
\newtheorem{prop}[thm]{Proposition}
\theoremstyle{remark} 
\newtheorem{rem}[thm]{Remark}
\theoremstyle{definition} 
\newtheorem{defn}[thm]{Definition}

\newcommand{\rl}{\mathbb{R}} 
\newcommand{\cx}{\mathbb{C}}

\newcommand{\go}{\omega} 
 
\newcommand{\gS}{\Sigma} 
\newcommand{\gp}{\varphi} 
\newcommand{\gn}{\nabla}

\newcommand{\ovl}[1]{\overline{#1}} 
 
\newcommand{\demi}{\frac{1}{2}} 
\newcommand{\gG}{\Gamma}

\title{Surfaces in $\bm{\mathbb{S}^3}$ and $\bm{\mathbb{H}^3}$ via Spinors}

\author{Bertrand Morel} 
\address{Bertrand Morel\newline
\indent Institut {\'E}lie Cartan\newline
\indent Universit\'e Henri Poincar\'e, Nancy I\newline
\indent B.P. 239\newline
\indent 54506 Vand\oe uvre-L{\`{e}}s-Nancy Cedex\newline
\indent France}
\email{morel@iecn.u-nancy.fr}



\begin{document}

\begin{abstract}
We generalize the spinorial characterization of isometric immersions of surfaces in
$\rl^3$ given in \cite{Fr} by T.~Friedrich to surfaces in $\mathbb{S}^3$ and
$\mathbb{H}^3$. The main argument is the interpretation of the
energy-momentum tensor associated with a special spinor field as a second
fundamental form. It turns out that such a characterization of isometric immersions
in terms of a special section of the spinor bundle also holds in the case
of hypersurfaces in the Euclidean $4$-space.
\end{abstract}
\maketitle

\section{Introduction}
It is well known that a description of a conformal immersion of an
arbitrary surface $M^2 \hookrightarrow \rl^3$ by a spinor field $\gp$ on
$M^2$ satisfying the inhomogenous Dirac equation 
\begin{equation}
  \label{eq:dirach}
 D\gp=H\gp, 
\end{equation}
(where $D$ stands for the Dirac operator and $H$ for the mean curvature of the surface), is possible. Recently,
many authors investigated such a description (see for example \cite{KuSch},\cite{Tai1}). 

In fact, it is clear that any oriented immersed surface $M^2
\hookrightarrow \rl^3$ inherits from $\rl^3$ a solution of Equation
(\ref{eq:dirach}), the surface $M$ being endowed with the induced metric and the
induced spin structure. Moreover, the solution has constant length. This
solution is obtained by the restriction to the surface of a parallel spinor field on
$\rl^3$. In
\cite{Fr}, T.~Friedrich clarifies the above-mentioned representation of surfaces
in $\rl^3$ in a geometrically invariant way by proving the following:

\begin{thm}[Friedrich \cite{Fr}]\label{Fried}
Let $(M^2,g)$ be an oriented, $2$-dimensional manifold and $H : M
\rightarrow \rl$ a smooth function. Then the following data are equivalent:
\begin{enumerate}
\item An isometric immersion $(\widetilde{M}^2,g)\rightarrow \rl^3$ of the
  universal covering $\widetilde{M}^2$ into the Euclidean space $\rl^3$ with
  mean curvature $H$.
\item A solution $\gp$ of the Dirac
  equation $$D\gp=H\gp\,,$$ with constant length $|\gp|\equiv 1$. 
\item A pair $(\gp,T)$ consisting of a symmetric endomorphism $T$ of the
  tangent bundle $TM$ such that
  $\mathrm{tr}(T)=H$ and a spinor field $\gp$ satisfying, for any $X\in\gG(TM)$, the equation $$\gn_X\gp+T(X)\cdot\gp=0\,.$$
\end{enumerate}
\end{thm}

In this paper, we prove the analoguous characterizations for surfaces in $\mathbb{S}^3$
and $\mathbb{H}^3$ (Theorems \ref{thms3} and \ref{thmh3}). They are obtained
by studying the equation of restrictions to a surface of real and imaginary Killing spinor fields (compare with \cite{Fr}). 

We note that the involved symmetric endomorphism $T$ is the energy-mo\-men\-tum tensor associated with the restricted
Killing spinor which describes the immersion. 

Finally, the case of the hypersurfaces of $\mathbb{R}^4$ is treated (Theorem \ref{thm53}).

\section{Restricting Killing spinor fields to a surface\label{secres}}

Let $N^3$ be a $3$-dimensional oriented Riemannian manifold, with a fixed
spin structure. Denote by $\gS N$ the spinor bundle associated with this spin
structure. If $M^2$ is an oriented surface isometrically immersed into
$N^3$, denote by $\nu$ its unit normal vector field. Then $M^2$ is endowed with a spin structure, canonically induced by
that of $N^3$. Denote by $\gS M$ the corresponding spinor bundle. The
following proposition is essential for what follows (see for example \cite{Bar1},\cite{BHMM},\cite{Mor},\cite{Trau1}):

\begin{prop}\label{ident}
There exists an identification of $\gS N_{|M}$ with $\gS M$, which after
restriction to $M$, sends
every spinor field $\psi \in \gG(\gS N)$ to the spinor field denoted by
$\psi^* \in \gG(\gS M)$. Moreover, if $\underset{N}{\cdot}$ (resp.
$\cdot$) stands for Clifford multiplication on $\gS N$ (resp. $\gS M$),
then one has
\begin{equation}
  \label{eq:clifmul}
  (X\underset{N}{\cdot}\nu\underset{N}{\cdot}\psi)^*=X\cdot\psi^*\;,
\end{equation}
for any vector field $X$ tangent to $M$.
\end{prop}

Another important formula is the well-known spinorial Gauss formula: if
$\gn^N$ and $\gn$ stand for the covariant derivatives on $\gG(\gS N)$ and
$\gG(\gS M)$  respectively, then, for all $X \in TM$ and $\psi \in
\gG(\gS N)$
\begin{equation}
  \label{eq:gauss}
  (\gn^N_X\psi)^*=\gn_X\psi^*+\demi h(X)\cdot\psi^*,
\end{equation}
where $h$ is the second fundamental form of the immersion $M\hookrightarrow
N$ viewed as a symmetric endomorphism of the tangent bundle of $M$.

Assume now that $N^3$ admits a non-trivial Killing spinor field of Killing constant
$\eta\in\cx$, i.e., a spinor field $\Phi\in\gG(\gS N)$ satisfying
\begin{equation}
  \label{eq:kill}
\gn^N_Y\Phi=\eta\, Y\underset{N}{\cdot}\Phi  
\end{equation}
for all vector field $Y$
on $N$. Recall that $\eta$ has to be real or pure imaginary and that
$\Phi$ never vanishes on $N$, as a non-trivial parallel section for a
modified connection (see \cite{BFGK},\cite{BHMM}). In what
follows, we will consider the model spaces, with
their standard metrics, $\rl^3$ with
$\eta=0$,  $\mathbb{S}^3$ with
$\eta=1/2$, and $\mathbb{H}^3$ with
$\eta=i/2$ which are characterized by the fact that they admit a
maximal number of linearly independant Killing spinor fields with constant
$\eta$.

Let $(e_1,e_2)$ be a positively oriented local orthonormal basis of $\gG(TM)$
such that $(e_1,e_2,\nu)$ is a positively oriented local orthonormal basis
of $\gG(TN)_{|M}$. Denote by
$$\go_3=-e_1\underset{N}{\cdot}e_2\underset{N}{\cdot}\nu$$ the complex volume
form on the complex Clifford bundle $\cx lN$ and $\go=e_1\cdot e_2$ the real volume
form on $\cx lM$. Recall that $\go_3$ acts by Clifford multiplication as
the identity on $\gS N$. Therefore, denoting $\gp:=\Phi^*$, formula
(\ref{eq:clifmul}) yields
$$(e_1\underset{N}{\cdot}\Phi)^*=(-e_1\underset{N}{\cdot}e_1\underset{N}{\cdot}e_2\underset{N}{\cdot}\nu\underset{N}{\cdot}\Phi)^*=e_2\cdot\gp=-e_1\cdot
\go\cdot\gp$$
$$(e_2\underset{N}{\cdot}\Phi)^*=(-e_2\underset{N}{\cdot}e_1\underset{N}{\cdot}e_2\underset{N}{\cdot}\nu\underset{N}{\cdot}\Phi)^*=-e_1\cdot\gp=-e_2\cdot\go\cdot\gp$$
and
$$(\nu\underset{N}{\cdot}\Phi)^*=(-\nu\underset{N}{\cdot}e_1\underset{N}{\cdot}e_2\underset{N}{\cdot}\nu\underset{N}{\cdot}\Phi)^*=\go\cdot\gp\;.$$

Then, these last relations with Equations (\ref{eq:gauss}) and
(\ref{eq:kill}) show that
\begin{equation}
  \label{rkill1}\forall X\in TM\,,\qquad\gn_X\gp+\demi h(X)\cdot\gp+\eta
  X\cdot\go\cdot\gp=0
\end{equation} 
Recall that the spinor bundle $\gS M$ splits into $$\gS
M=\gS^+M\oplus\gS^-M$$ where $\gS^\pm M$ is the $\pm 1$-eigenspace for the
action of the complex volume forme $\go_2=i\,\go$. Under this
decomposition, we will denote $\gp=\gp^+ +\gp^-$, and define
$\ovl{\gp}:=\gp^+-\gp^-$. Therefore Equation (\ref{rkill1}) is equivalent to 
$$\gn_X\gp+\demi h(X)\cdot\gp-i\eta X\cdot\ovl{\gp}=0\;.$$

The ambient spinor bundle $\gS N$ can be endowed  with a Hermitian inner product $(\;,)_N$ for which Clifford
multiplication by any vector tangent to $N$ is skew-symmetric. This product
induces another Hermitian inner product on $\gS M$, denoted by $(\;,)$ making the
identification of Proposition \ref{ident} an isometry. Now, relation
(\ref{eq:clifmul}) shows that Clifford
multiplication by any vector tangent to $M$ is skew-symmetric with respect
to $(\;,)$. 

\begin{prop}\label{lengthphi}
If $\eta\in\rl$, then $\gp$ has constant length. If $\eta\in
i\rl^*$, then for all vector $X$ tangent to $M$, $$X|\gp|^2=2\Re(i\eta X\cdot\ovl{\gp},\gp)\;.$$
\end{prop}

\begin{proof} Since Clifford
multiplication by any vector tangent to $M$ is skew-symmetric with respect
to $(\;,)$, we have $\Re(Y\cdot\gp,\gp)=0$ for all $Y\in TM$. Taking this
fact into account and computing $$X|\gp|^2=2\Re(\gn_X\gp,\gp)$$ with the
help of formula (\ref{rkill1}), completes the proof.
\end{proof}

Recalling that the Dirac operator $D$ is defined on $\gG(\gS M)$ by 
$$D=e_1\cdot\gn_{e_1}+e_2\cdot\gn_{e_2}\;,$$
we compute directly that $$D\gp=H\gp+2\eta\,\go\cdot\gp=H\gp-2i\eta\ovl{\gp}$$
where $H$ is the mean curvature of the immersion $M\hookrightarrow
N$.  It is well known that the
action of the Dirac operator satisfies $(D\gp)^\pm=D\gp^\mp$ (see \cite{LM},\cite{BHMM}). Therefore, we note
that 
\begin{equation}
  \label{eq:dgp}
D(\gp^\pm)=(H\pm 2i\eta)\gp^\mp\;.  
\end{equation}

We have as in \cite{Fr} the following

\begin{prop}
Let $M^2$ be a minimal surface in $N^3$. Then the restriction of any
Killing spinor $\Phi$ with constant $\eta$ on $N^3$ restricts to an eigenspinor $\gp^\star$ on the
surface $M$: $$D\gp^\star=2\eta\gp^\star$$ Moreover, if $\eta$ is
real, then $\gp^\star$ has constant length.
\end{prop}

\begin{proof}
Since $H=0$, we have $$D(\gp^\pm)=\pm 2i\eta\gp^\mp\;.$$ Therefore, it
suffices to define $\gp^\star=\gp^++i\gp^-$. 
\end{proof}

\section{Solutions of the restricted Killing spinor equation}

Let $(M^2,g)$ be an oriented, 2-dimensional Riemannian manifold with a spin
structure. We endow the spinor bundle $\gS M$ with a Hermitian inner product $(\;,)$
for which Clifford
multiplication by any vector tangent to $M$ is skew-symmetric.

We study now some properties of a given solution $\gp\in\gG(\gS M)$ of the following equation
\begin{equation}
  \label{eq:restkill}
 \gn_X\gp+T(X)\cdot\gp-i\eta X\cdot\ovl{\gp}=0\; , 
\end{equation} 
or equivalently
\begin{equation}
  \label{eq:restkill2}
 \gn_X\gp+T(X)\cdot\gp+\eta X\cdot\go\cdot\gp=0\; , 
\end{equation} 
where $T$ stand for a symmetric endomorphism of the tangent bundle of $M$, and $\eta
\in \rl\cup i\rl$. 

In view of the preceding section and for reasons which will become clearer later, we will call this
equation the {\it restricted Killing spinor equation}\,.
The following proposition shows the role of solutions of the restricted Killing spinor equation in the theory of surfaces in 
$\rl^3$, $\mathbb{S}^3$ and $\mathbb{H}^3$. In fact, we see that the integrability conditions for such sections of the spinor bundle are 
exactly the Gau{\ss} and Codazzi-Mainardi equations. 

In the following, $(e_1,e_2)$ denotes a positively oriented local orthonormal basis of
$\gG(TM)$.

\begin{prop}\label{gaussmain}
Assume that $(M^2,g)$ admits a non trivial solution of Equation
(\ref{eq:restkill}) and let $S=2T$, then 
$$(\gn_XS)(Y)=(\gn_YS)(X)\qquad\text{(Codazzi-Mainardi Equation),}$$and
$$R_{1212}-\det(S)=4\eta^2\qquad\text{(Gau{\ss} Equation),}$$
where $R_{1212}=g(R(e_1,e_2)\,e_2,e_1)$, and $R$ is the Riemann tensor of $M$.\end{prop} 

\begin{proof}
Let $\gp$ a non-trivial solution of (\ref{eq:restkill}). We compute the action of the spinorial curvature tensor
$\mathcal{R}$ on $\gp$ defined
for all $X,Y \in TM$ by 
$$\mathcal{R}(X,Y)\gp=\gn_X\gn_Y\gp-\gn_Y\gn_X\gp-\gn_{[X,Y]}\gp\;.$$
Since it is skew-symmetric and $\dim M=2$, with the
help of formula (\ref{eq:restkill2}), we only compute
\begin{eqnarray*}
\gn_{e_1}\gn_{e_2}\gp&=&\gn_{e_1}(-T(e_2)\cdot\gp-\eta e_2\cdot\go\cdot\gp)\\
&=&\gn_{e_1}(-T(e_2)\cdot\gp-\eta e_1\cdot\gp)\\
&=&-\gn_{e_1}T(e_2)\cdot\gp-T(e_2)\cdot\gn_{e_1}\gp-\eta\gn_{e_1} e_1\cdot\gp-\eta
e_1\cdot\gn_{e_1}\gp\\
&=&-\gn_{e_1}T(e_2)\cdot\gp+T(e_2)\cdot T(e_1)\cdot\gp-\eta T(e_2)\cdot e_2\cdot\gp\\&&-\eta\gn_{e_1} e_1\cdot\gp+\eta
e_1\cdot T(e_1)\cdot\gp-\eta^2 e_1\cdot e_2\cdot \gp
\end{eqnarray*} as well as
\begin{eqnarray*}\gn_{e_2}\gn_{e_1}\gp&=&-\gn_{e_2}T(e_1)\cdot\gp+T(e_1)\cdot T(e_2)\cdot\gp+\eta T(e_1)\cdot e_1\cdot\gp\\&&+\eta\gn_{e_2} e_2\cdot\gp-\eta
e_2\cdot T(e_2)\cdot\gp+\eta^2 e_1\cdot e_2\cdot\gp\;.\end{eqnarray*}
So, taking into account that $[e_1,e_2]=\gn_{e_1}e_2-\gn_{e_2}e_1$, a
straightforward computation gives
\begin{eqnarray}
 \mathcal{R} (e_1,e_2)\gp&=&\Big( (\gn_{e_2}T)(e_1)-(\gn_{e_1}T)(e_2)\Big )\cdot\gp\nonumber\\
&&-\Big(T(e_1)\cdot T(e_2)-T(e_2)\cdot T(e_1)\Big)\cdot\gp\nonumber\\
&&-2\eta^2 e_1\cdot e_2\cdot \gp\label{presqr} \end{eqnarray}
On the other hand, it is well known that this spinorial curvature tensor
corresponds to the Riemann tensor $R$ of $M$ via the relation
\begin{equation}\label{RR}\mathcal{R}(e_1,e_2)\gp=-\demi R_{1212}e_1\cdot e_2\cdot\gp\;.\end{equation}
Now, it is easy to see that $$T(e_1)\cdot T(e_2)-T(e_2)\cdot
T(e_1)=2\det(T)e_1\cdot e_2$$ and therefore, if we put $S=2T$ and define
the {\it function} $$G:=R_{1212}-\det(S)-4\eta^2$$ and the {\it vector
  field} $$C:=(\gn_{e_1}S)(e_2)-(\gn_{e_2}S)(e_1),$$ Equations (\ref{presqr})
and (\ref{RR}) yield $$C\cdot\gp=G  e_1\cdot e_2\cdot \gp\;.$$ Note that $e_1\cdot
e_2\cdot \gp=-i\,\ovl{\gp}$, hence $$C\cdot\gp^\pm=\pm iG\gp^\mp\;.$$
Applying two times this relation, it
suffices to note that
$$||C||^2\gp^\pm=-G^2\gp^\pm\;,$$ and so $C=0$ and $G=0$.
\end{proof}

Note that up to rescaling, we
can take $\eta=0$, $1/2$, or $i/2$. The case $\eta=0$ is treated in
\cite{Fr} and is the starting point of the proof of Theorem \ref{Fried}. We
will discuss the cases $\eta=1/2$ and $\eta=i/2$ separately.
We begin by

\begin{lem}\label{defE}
Let $\gp$ be a non trivial solution of the restricted Killing spinor equation
(\ref{eq:restkill}). Then 
\begin{itemize}
\item if $\eta=1/2$, $\gp$ has constant norm and the symmetric 
endomorphism $T$, viewed as a covariant symmetric $2$-tensor, is
given by 
$$T(X,Y)=\demi\Re(X\cdot\gn_Y\gp+Y\cdot\gn_X\gp,\gp/|\gp|^2)$$
\item if $\eta=i/2$, $\gp$ satisfies $X|\gp|^2=-\Re(X\cdot\ovl{\gp},\gp)$ and one has
\end{itemize} 
$$T(X,Y)|\gp|^2=\demi\Re(X\cdot\gn_Y\gp+Y\cdot\gn_X\gp,\gp)+\demi\Big(|\gp^-|^2-|\gp^+|^2\Big)g(X,Y)$$
\end{lem}

\begin{proof}
The first claim of each case is proved in Proposition \ref{lengthphi}. Let $T_{jk}=g(T(e_j),e_k)$, then, for $j=1,2$, $$\gn_{e_j}\gp=-\sum_{k=1}^2
T_{jk}e_k\cdot\gp+i\eta e_j\cdot\ovl{\gp}\;.$$
Taking Clifford multiplication by $e_l$ and the scalar product with $\gp$,
we get
$$\Re(e_l\cdot\gn_{e_j}\gp,\gp)=-\sum_{k=1}^2 T_{jk}\Re(e_l\cdot
e_k\cdot\gp,\gp)+\Re(i\eta e_l\cdot e_j\cdot\ovl{\gp},\gp)\;.$$
Since $\Re(e_l\cdot
e_k\cdot\gp,\gp)=-\delta_{lk}|\gp|^2$, it follows, by symmetry of $T$
$$\Re(e_l\cdot\gn_{e_j}\gp+e_j\cdot\gn_{e_l}\gp,\gp)=2T_{lj}|\gp|^2-2\Re(i\eta\ovl{\gp},\gp)\delta_{lj}\;.$$
This completes the proof by taking $\eta=1/2$ or  $\eta=i/2$.
\end{proof}

Now, we prove that the necessary conditions on a spinor field $\psi\in\gG(\gS
M)$ obtained in the previous section
(i.e. Proposition \ref{lengthphi} and Equation (\ref{eq:dgp}))
are enough to prove that $\psi$ is a solution of the restricted Killing
spinor equation.

\vspace{3mm}
\noindent{\bf The case $\bm{\eta=1/2}$:} Consider a non-trivial
spinor field $\psi$ of constant length, which satisfies $D\psi^\pm=(H\pm
i)\psi^\mp$. Define the following $2$-tensors on $(M^2,g)$
$$T^\pm(X,Y)=\Re(\gn_X\psi^\pm,Y\cdot\psi^\mp)\;.$$ First note that 
\begin{equation}
  \label{eq:trE}
 \mathrm{tr}T^\pm=-\Re(D\psi^\pm,\psi^\mp)=-\Re((H\pm
i)\psi^\mp,\psi^\mp)=-H|\psi^\mp|^2\;. 
\end{equation}
We have the following relations
\begin{eqnarray}
 T^\pm(e_1,e_2)&=&
 \Re(\gn_{e_1}\psi^\pm,e_2\cdot\psi^\mp)=\Re(e_1\cdot\gn_{e_1}\psi^\pm,e_1\cdot e_2\cdot\psi^\mp)\nonumber\\
&=&\Re(D\psi^\pm,e_1\cdot e_2\cdot\psi^\mp)-\Re(e_2\cdot\gn_{e_2}\psi^\pm,e_1\cdot e_2\cdot\psi^\mp)\nonumber\\
&=&\Re((H\pm i)\psi^\mp,e_1\cdot
e_2\cdot\psi^\mp)+\Re(\gn_{e_2}\psi^\pm,e_1\cdot\psi^\mp)\nonumber\\
&=&|\psi^\mp|^2+ T^\pm(e_2,e_1)\;.\label{eq:passym}\end{eqnarray}

\begin{lem}\label{relE}
The $2$-tensors $T^\pm$ are related by the equation
$$|\psi^+|^2T^+=|\psi^-|^2T^-$$
\end{lem}

\begin{proof}
This relation is trivial at any point $p \in M$ where $|\psi^+|^2$ or
$|\psi^-|^2$ vanishes. Therefore we can assume in the following that both spinors $ \psi^+$
and $\psi^-$ are not zero in the neighbourhood of a point in $M$. 

With respect to the scalar product $\Re(\,, )$, the spinors
$$e_1\cdot\frac{\psi^-}{|\psi^-|}\quad\text{  and
  }\quad e_2\cdot\frac{\psi^-}{|\psi^-|}$$ form a local orthonormal basis of $\gG(\gS^+
M)$. Hence, in this basis, we can write
\begin{eqnarray*}
  \gn_X\psi^+&=&\Re(\gn_X\psi^+,e_1\cdot\frac{\psi^-}{|\psi^-|}
  )\,e_1\cdot\frac{\psi^-}{|\psi^-|}+\Re(\gn_X\psi^+,e_2\cdot\frac{\psi^-}{|\psi^-|} )\,e_2\cdot\frac{\psi^-}{|\psi^-|}\\
&=&\frac{T^+(X)}{|\psi^-|^2}\cdot\psi^-
\end{eqnarray*} where the vector field $T^+(X)$ is defined by
$$g(T^+(X),Y)=T^+(X,Y)\,,\qquad\forall Y\in TM\;.$$ In the same manner, we can show
that $$\gn_X\psi^-=\frac{T^-(X)}{|\psi^+|^2}\cdot\psi^+\;.$$
Since $\psi$ has constant length, for all vector $X$ tangent to $M$, we have 
\begin{eqnarray}
  0&=&X|\psi|^2=X(|\psi^+|^2+|\psi^-|^2)\nonumber\\
&=&2\Re(\gn_X\psi^+,\psi^+)+2\Re(\gn_X\psi^-,\psi^-)\nonumber\\
&=&2\Re(W(X)\cdot\psi^-,\psi^+)\label{rank}\end{eqnarray} with $$W(X)=\frac{T^+(X)}{|\psi^-|^2}-\frac{T^-(X)}{|\psi^+|^2}\;.$$
To conclude, it suffices to note that Equations
(\ref{eq:trE}) and (\ref{eq:passym}) imply $W$ is traceless and
symmetric, and that Equation (\ref{rank}) implies that $W$ has rank less
or equal to $1$. This obviously implies $W=0$.
\end{proof} 

\begin{prop}
Assume that there exists on $(M^2,g)$ a non-trivial solution $\psi$ of the
equation $D\psi=H\psi-i\ovl{\psi}$ with constant length. Then such a solution satisfies the restricted Killing spinor equation with $\eta=1/2$.
\end{prop}

\begin{proof}
Let $F:=T^++T^-$. Lemma \ref{relE} and the begining of its proof imply
$$\frac{F}{|\psi|^2}=\frac{T^+}{|\psi^-|^2}=\frac{T^-}{|\psi^+|^2}\;.$$ Hence
$F/|\psi|^2$ is well defined on the whole surface $M$, and \begin{equation}\label{preqrsk}\gn_X\psi=\gn_X\psi^++\gn_X\psi^-=\frac{F(X)}{|\psi|^2}\cdot\psi\end{equation}
 where the vector field $F(X)$ is defined by
$g(F(X),Y)=F(X,Y)$, $\forall Y\in TM.$ Note that by Equation
(\ref{eq:passym}), the $2$-tensor $F$ is not symmetric. Define now the symmetric $2$-tensor
$$T(X,Y)=-\frac{1}{2|\psi|^2}\left(F(X,Y)+F(Y,X)\right)\;.$$ Observe that $T$ is
defined as in Lemma \ref{defE}. It is straigthforward to
show that
$$T(e_1,e_1)=-F(e_1,e_1)/|\psi|^2\quad,\quad T(e_2,e_2)=-F(e_2,e_2)/|\psi|^2\;,$$
$$T(e_1,e_2)=-F(e_1,e_2)/|\psi|^2+\demi\quad\text{and}\quad T(e_2,e_1)=-F(e_2,e_1)/|\psi|^2-\demi$$ once more by Equation
(\ref{eq:passym}). Taking into account these last relations in Equation
(\ref{preqrsk}), we conclude $$\gn_X\psi=-T(X)\cdot\psi-\demi X\cdot\go\cdot\psi\;.$$
\end{proof}

\noindent{\bf The case $\bm{\eta=i/2}$:}

\begin{prop}
Assume that there exists on $(M^2,g)$ a nowhere vanishing solution $\psi$ of the equation
$D\psi=H\psi+\ovl{\psi}$. Then, if this solution satisfies $$X|\psi|^2=-\Re(X\cdot\ovl{\psi},\psi)\;,\qquad \forall X\in\gG(TM),$$
then it is solution of the restricted Killing spinor equation with $\eta=i/2$.
\end{prop}

\begin{proof}
 Defining the $2$-tensors $T^\pm$ as in the previous case,
we get 
\begin{equation}
  \label{eq:traceEi}
\mathrm{tr}T^\pm=-(H\mp 1)|\psi^\mp|^2\;,
\end{equation}
 and 
 \begin{equation}
   \label{eq:Esym}
 T^\pm(e_1,e_2)=T^\pm(e_2,e_1)\;.  
 \end{equation}
First note that
$$-\Re(X\cdot\ovl{\psi},\psi)=-\Re(X\cdot\psi^+,\psi^-)+\Re(X\cdot\psi^-,\psi^+)=2\Re(X\cdot\psi^-,\psi^+)\;.$$
Therefore, following the proof of Lemma \ref{relE}, we get
\begin{equation}
  \label{eq:aeqI}
\Re(X\cdot\psi^-,\psi^+)=\Re(W(X)\cdot\psi^-,\psi^+)
\end{equation}
 with $$W(X)=\frac{T^+(X)}{|\psi^-|^2}-\frac{T^-(X)}{|\psi^+|^2}\;.$$
As in the previous case, Equations (\ref{eq:traceEi}), (\ref{eq:Esym}) and
(\ref{eq:aeqI}) imply that $W-\mathrm{Id}_{TM}$ is a symmetric, traceless
endomorphism of rank not greater than $1$, hence $W=\mathrm{Id}_{TM}$ and we have
the relation $$|\psi^+|^2T^+-|\psi^-|^2 T^- =|\psi^+|^2|\psi^-|^2g\;.$$
Therefore, if we define the symmetric $2$-tensor
$F=T^++T^-+\demi(|\psi^+|^2-|\psi^-|^2)g$, we have on the whole surface $M$
$$\frac{F}{|\psi|^2}=\frac{T^++T^-+(|\psi^+|^2-|\psi^-|^2)g}{|\psi^+|^2+|\psi^-|^2}=\frac{T^-}{|\psi^+|^2}+\demi
g=\frac{T^+}{|\psi^-|^2}-\demi g\;.$$
On the other hand, we get
$$\gn_X\psi=\gn_X\psi^++\gn_X\psi^-=\frac{T^+(X)}{|\psi^-|^2}\cdot\psi^-+\frac{T^-(X)}{|\psi^+|^2}\cdot\psi^+\;.$$
These two last equations imply
$$\gn_X\psi=\frac{F(X)}{|\psi|^2}\cdot(\psi^++\psi^-)+\demi X\cdot\psi^--\demi
X\cdot\psi^+\;,$$ which is equivalent to $$\gn_X\psi=-T(X)\cdot\psi-\demi
X\cdot\ovl{\psi}\;.$$ Naturally, we put $T=-\frac{F}{|\psi|^2}$ and note
that $T$ is defined as in Lemma \ref{defE}.
\end{proof}

\section{Surfaces in $\mathbb{S}^3$ or $\mathbb{H}^3$}

We are now able to generalize Theorem \ref{Fried} to surfaces in $\mathbb{S}^3$ or
$\mathbb{H}^3$. In section 2, we saw that an oriented, immersed surface
$M^2\hookrightarrow \mathbb{S}^3$ (resp. $\mathbb{H}^3$) inherits an induced metric $g$, a
spin structure, and a solution $\gp$ of 
\begin{equation}\label{dirdir}
D\gp=H\gp-i\ovl{\gp}\qquad
\text{(resp. }D\gp=H\gp+\ovl{\gp}\;\text{)}
\end{equation} 
with constant length
(resp. with $X|\gp|^2=-\Re(X\cdot\ovl{\gp},\gp)$ for all vector $X$ tangent
to $M$). This spinor field $\gp$ on $M^2$ is the restriction of a real
(resp. imaginary) Killing spinor field in $\mathbb{S}^3$ (resp. $\mathbb{H}^3$).
Section 3 shows that at least locally the converse is true. Assume that there exists a
solution of Equation \eqref{dirdir} on an oriented, $2$-dimensional Riemannian
manifold $(M^2,g)$ endowed with a spin structure, for a given function
$H:M\rightarrow \rl$. Then this solution
satisfies the restricted Killing spinor equation with a well defined
endomorphism $T:TM \rightarrow TM$ with $\mathrm{tr}T=H$. Moreover, there
exists an isometric immersion $(M^2,g)\hookrightarrow \mathbb{S}^3$ (resp. $\mathbb{H}^3$)
with second fundamental form $S=2T$.

\begin{thm}\label{thms3}
Let $(M^2,g)$ be an oriented, $2$-dimensional manifold and $H : M
\rightarrow \rl$ a smooth function. Then the following data are equivalent:
\begin{enumerate}
\item An isometric immersion $(\tilde{M}^2,g)\rightarrow \mathbb{S}^3$ of the
  universal covering $\tilde{M}^2$ into the $3$-dimensional round sphere $\mathbb{S}^3$ with
  mean curvature $H$.
\item A solution $\gp$ of the Dirac
  equation $$D\gp=H\gp-i\ovl{\gp}$$ with constant length. 
\item A pair $(\gp,T)$ consisting of a symmetric endomorphism $T$ such that
  $\mathrm{tr}(T)=H$ and a spinor field $\gp$ satisfying the equation $$\gn_X\gp+T(X)\cdot\gp-\frac{i}{2}X\cdot\ovl{\gp}=0\,.$$
  \end{enumerate}
\end{thm}
\begin{thm}\label{thmh3}
Let $(M^2,g)$ be an oriented, $2$-dimensional manifold and $H : M
\rightarrow \rl$ a smooth function. Then the following data are equivalent:
\begin{enumerate}
\item An isometric immersion $(\tilde{M}^2,g)\rightarrow \mathbb{H}^3$ of the
  universal covering $\tilde{M}^2$ into the $3$-dimensional hyperbolic space $\mathbb{H}^3$ with
  mean curvature $H$.
\item A nowhere vanishing solution $\gp$ of the Dirac
  equation $$D\gp=H\gp+\ovl{\gp}$$ satisfying $$X|\gp|^2=-\Re(X\cdot\ovl{\gp},\gp)\quad \forall X\in\gG(TM).$$
\item A pair $(\gp,T)$ consisting of a symmetric endomorphism $T$ such that
  $\mathrm{tr}(T)=H$ and a spinor field $\gp$ satisfying the equation $$\gn_X\gp+T(X)\cdot\gp+\frac{1}{2}X\cdot\ovl{\gp}=0\quad \forall X\in\gG(TM).\,.$$
\end{enumerate}
\end{thm}

\begin{rem}
It has been pointed out to us that the case of surfaces in $\mathbb{S}^3$ has
already been treated by Leonard Voss ({\it Diplomarbeit,
  Humboldt-Universit\"at zu Berlin}, unpublished). 
\end{rem}

\section{Hypersurfaces in $\rl^4$}

We conclude by giving a characterization of hypersurfaces in the Euclidean $4$-space in terms of a special section of
the intrinsic spinor bundle of the hypersurface, in a very similar way to
that of Theorem \ref{Fried}. 

Let $M^3$ be an oriented hypersurface isometrically immersed into
$\mathbb{R}^4$, denote by $\nu$ its unit normal vector field. Then $M^3$ is endowed with a spin structure, canonically induced by
that of $\mathbb{R}^4$. Denote by $\gS M$ the corresponding spinor
bundle and $\gS^+\mathbb{R}^4$ the bundle of positive spinors in
$\mathbb{R}^4$. We then have the anologous result of Proposition \ref{ident}:

\begin{prop}\label{ident2}
There exists an identification of $\gS^+\mathbb{R}^4$ with $\gS M$, which after
restriction to $M$, sends
every spinor field $\psi \in \gG(\gS^+\mathbb{R}^4)$ to the spinor field denoted by
$\psi^* \in \gG(\gS M)$. Moreover, if $\underset{\mathbb{R}^4}{\cdot}$ (resp.
$\cdot$) stands for Clifford multiplication on $\gS^+\mathbb{R}^4$ (resp. $\gS M$),
then one has
\begin{equation}
  (X\underset{\mathbb{R}^4}{\cdot}\nu\underset{\mathbb{R}^4}{\cdot}\psi)^*=X\cdot\psi^*\;,
\end{equation}
for any vector field $X$ tangent to $M$.
\end{prop}

Recall the following definition

\begin{defn} A symmetric $2$-tensor $T\in S^2(M)$ is called a \emph{Codazzi
    tensor} if it satisfies the Codazzi-Mainardi equation, i.e. 
$$(\gn_XT)(Y)=(\gn_YT)(X)\qquad\forall X,Y\in\gG(TM)\;,$$ 
($T$ being viewed in this formula via the metric $g$ as a symmetric endomorphism of the
tangent bundle).
\end{defn}

We now prove the following

\begin{thm}\label{thm53}
Let $(M^3,g)$ be an oriented, $3$-dimensional Riemannian manifold. Then the following data are equivalent:
\begin{enumerate}
\item An isometric immersion $(\widetilde{M}^3,g)\rightarrow \rl^4$ of the
  universal covering $\widetilde{M}^3$ into the Euclidean space $\rl^4$
  with second fundamental form $h$. 
\item A pair $(\gp,T)$ consisting of a Codazzi tensor $T$ such that
  $2T=h$ and a non trivial spinor field $\gp$ satisfying, for any $X\in\gG(TM)$, the equation $$\gn_X\gp+T(X)\cdot\gp=0\,.$$
\end{enumerate}
\end{thm}

\begin{proof}
Let $(M^3,g)$ be an oriented hypersurface isometrically immersed into
$\rl^4$ with second fundamental form $h$. Let $\psi$ be any parallel
positive spinor
field on $\rl^4$. Denote by $\gp:=\psi^*\in\gG(\gS M)$ the restriction of
$\psi$ given by Proposition \ref{ident2}. Then Gau{\ss} formula (\ref{eq:gauss}) yields 
$$\gn_X\gp+\demi h(X)\cdot\gp=0\;.$$
Since $h$ is a second fundamental form, it is clear that $T=\demi h$ is a
Codazzi tensor and that $(\gp,T)$ give the desired pair.

Conversely, if $(M^3,g)$ is an oriented, $3$-dimensional Riemannian
manifold admitting such a pair $(\gp,T)$, then obviously Codazzi-Mainardi
equation holds for $h=2T$. 

Therefore, the action of the spinorial curvature tensor on the spinor $\gp$ is given by 
\begin{equation}\label{eqaction}
  \mathcal{R}(X,Y)\gp = \Big ( T(Y)\cdot T(X)-T(X)\cdot T(Y)\Big
  )\cdot\gp
\end{equation}

Let $(e_1,e_2,e_3)$ be a positively oriented local orthonormal basis of
$\gG(TM)$. Then Equation \eqref{eqaction} yields 
$$\sum_{k\neq l}\Big (\mathcal{R}_{ijkl} + 4T_{il}T_{jk}-4T_{ik}T_{jl}\Big
  )e_k \cdot e_l \cdot\gp=0$$ which imply in dimension $3$ that each
  componant $$\mathcal{R}_{ijkl} + 4T_{il}T_{jk}-4T_{ik}T_{jl}$$ is zero,
  since for $1\leq k < l\leq 3$ and $1\leq k' < l'\leq 3$, $$\Re (e_k \cdot
  e_l \cdot\gp,e_{k'} \cdot e_{l'}
  \cdot\gp)=\pm\delta_{kk'}\delta_{ll'}|\gp|^2\;.$$ Therefore $h=2T$
  satisfies the Gau\ss{} equation.
\end{proof}

\begin{rem}
Let $(\gp,T)$ be a pair as in Theorem \ref{thm53} ($2$). Then necessarily
the Codazzi tensor $T$ has to be defined as the energy-momentum tensor
associated with the spinor field $\gp$ (see for example \cite{Hi3},
\cite{FK1} or \cite{MoEM}). Such a special spinor field is then
called a \emph{Codazzi Energy-Momentum spinor}, and generalizes the notion
of Killing spinors (see \cite{MoCEM} for a study of these particular spinor
fields).
\end{rem}

\providecommand{\bysame}{\leavevmode\hbox to3em{\hrulefill}\thinspace}

\end{document}